\documentclass[12pt,leqno]{amsart}

\usepackage{color}

\newtheorem{theorem}{Theorem}
\newtheorem{proposition}[theorem]{Proposition}
\newtheorem{corollary}[theorem]{Corollary}
\newtheorem{remark}{Remark}
\newtheorem{lemma}[theorem]{Lemma}
\newfont{\bb}{msbm10 at 12pt}

\def\pf{\noindent{\textit {Proof.} }}

\def\R{\hbox{\bb R}}

\def\S{\hbox{\bb S}}

\def\H{\hbox{\bb H}}

\def\Ric{\hbox{Ric}}

\newcommand{\bal}{\begin{align}}      \newcommand{\eal}{\end{align}}
\newcommand{\ba}{\begin{array}}      \newcommand{\ea}{\end{array}}
\newcommand{\bc}{\begin{center}}     \newcommand{\ec}{\end{center}}
\newcommand{\be}{\begin{enumerate}}  \newcommand{\ee}{\end{enumerate}}
\newcommand{\beq}{\begin{eqnarray}}  \newcommand{\eeq}{\end{eqnarray}}
\newcommand{\beQ}{\begin{eqnarray*}} \newcommand{\eeQ}{\end{eqnarray*}}
\newcommand{\bi}{\begin{itemize}}    \newcommand{\ei}{\end{itemize}}
\newcommand{\bt}{\begin{tabular}}    \newcommand{\et}{\end{tabular}}
\newcommand{\bdm}{\begin{displaymath}} \newcommand{\edm}{\end{displaymath}}

\def\qed{\hfill{q.e.d.}\smallskip\smallskip}

\begin{document}

\title[Uniqueness of the AdS spacetime]{Uniqueness of the AdS spacetime among static vacua with prescribed
null infinity}

\author{Oussama Hijazi}
\address[Hijazi]{Institut {\'E}lie Cartan\\
Universit{\'e} de Lorraine, Nancy I\\
B.P. 239\\
54506 Vand\oe uvre-L{\`e}s-Nancy Cedex, France}
\email{Oussama.Hijazi@univ-lorraine.fr}

\author{Sebasti{\'a}n Montiel}
\address[Montiel]{Departamento de Geometr{\'\i}a y Topolog{\'\i}a\\
Universidad de Granada\\
18071 Granada \\
Spain}
\email{smontiel@ugr.es}

\begin{abstract}
We prove that an $(n+1)$-dimensional spin static vacuum with negative cosmological constant whose null infinity has a boundary admitting a non-trivial Killing spinor field is the AdS spacetime. As a consequence, we generalize  
previous uniqueness results by X. Wang \cite{Wa2}  and by Chru{\'s}ciel-Herzlich \cite{CH} and introduce,
for this class of spin static vacua, some Lorentzian manifolds which are prohibited as null infinities. 
\end{abstract}



\thanks{The second author was partially 
supported by a Spanish MEC-FEDER grants Nos. MTM2007-61775 and MTM2011-22547}

\date{November 23, 2012}        

\maketitle \pagenumbering{arabic}
 
\section{Introduction}
An $(n+1)$-dimensional {\em vacuum} spacetime with  cosmological constant $\Lambda< 0$ is a Lorentzian manifold $({\mathcal V},
g_{\mathcal V})$ satisfying the Einstein equation ${\Ric}_{\mathcal V}=\Lambda\,{g}_{\mathcal V}$. The vacuum is said to be
{\em static} when there exists a global time function $t$ with respect to which the spacelike slices $t=t_0$ are isometric
to each other, that is, space remains unchanged when time passes. The fact that $({\mathcal V},
g_{\mathcal V})$ is static imposes a precise kind of topologies and metrics, that is, imposes
the following form for the corresponding spacetime: 
\begin{equation}\label{static}
{\mathcal V}=\R\times M,\qquad g_{\mathcal V}=-V^2\,dt^2+g,
\end{equation}
where $(M,g)$ is an $n$-dimensional connected Riemannian manifold standing for the 
unchanging slices of constant time and $V\in C^\infty (M)$ is a positive smooth function
 on $M$. Without loss of generality, we can assume that $\Lambda = -n$. Then,  the normalized 
 vacuum Einstein equation 
 $\Ric_{\mathcal V}=-n g_{\mathcal V}$,
can be translated to the following two conditions on $(M,g)$ and $V$: 
\begin{equation}\label{vacuum}
\hbox{\rm Ric}_g+ng=\frac{\nabla^2V}{V},\qquad \Delta V=nV,
\end{equation}
where $\Ric_g$ is the Ricci tensor, $\nabla^2$ is the Hessian operator and $\Delta=
\hbox{\rm trace\,}\nabla^2$
is the Laplacian of the Riemannian manifold $(M,g)$. Taking traces in
the first of these two equations and taking into account the second one, we conclude
immediately that 
\begin{equation}\label{scalar}
R_g=-n(n-1),
\end{equation} 
where $R_g$ is the scalar curvature of $(M,g)$. Since 
the Riemannian manifold $(M,g)$ and the function $V$ completely determine the 
vacuum spacetime $({\mathcal V},g_{\mathcal V})$, it is also usual to call the triple
$(M,g,V)$ a {\em static vacuum}. \\ 

The paradigmatic example of a static vacuum is given by choosing $(M,g)=(\H^n,g_{\mathbb H})$ the hyperbolic space of constant sectional curvature $-1$, which can be viewed as
the upper sheet of a hyperquadric in the Minkowski spacetime $(\R_1^{n+1},\langle\;,\;\rangle)$, namely,$$
\H^n=\{p\in\R_1^{n+1}\,|\, |p|^2=-1, p_0>0\}.$$ 
This is a spacelike hypersurface and the  metric $g_{\Bbb H}$ induced on
$\H^n$ from $\langle\;,\;\rangle$ is a Riemannian metric making $(\H^n,g_{\Bbb H})$
a  space form of sectional curvature $-1$. In particular, $\Ric_{g_{\mathbb H}}=-(n-1)g_{\Bbb H}$. In this representation of the hyperbolic space, it is
easy to see that, for each fixed $a\in\R_1^{n+1}$, the {\em height} function $h_a:\H^n\rightarrow\R$ given by$$
h_a(p)=\langle p,a\rangle,\quad \forall p\in\H^n$$
satisfies the Obata type equation (see \cite{K})$$
\nabla^2 h_a=h_a\, g_{\Bbb H}.$$
As a consequence, $\Delta h_a=nh_a$. Moreover, when $a$ is chosen to be  lightlike or timelike vector with $a_0<0$, we have $h_a>0$.
Then the triple $(\H^n, g_{\Bbb H},h_a)$ is a static vacuum for such non-spacelike $a\in\R_1^{n+1}$.
Indeed, the corresponding spacetime $(AdS, g_{AdS})= (\R\times\H^n, -h_a^2\,dt^2+g_{\Bbb H})$ is nothing but the
Anti de Sitter (AdS) vacuum. \\ 

However, there is a subtle difference between the case where 
$a$ is timelike and the case where $a$ is lightlike. In the former, the level hypersurfaces of $h_a$ are umbilical round hyperspheres in $\H^n$ allowing to write the hyperbolic metric $g_{\mathbb H}$ as a warped product metric and describe the function $h_a$ as follows:
$$
g_{\mathbb H}=dr^2+\sinh^2 r\,\langle\;,\;\rangle_{{\mathbb S}^{n-1}},\quad h_a=\cosh r,
\qquad r\in[0,+\infty[.$$ 
Consequently, $$
g_{AdS}=-\cosh^2r\,dt^2+dr^2+\sinh^2r\,\langle\;,\;\rangle_{{\mathbb S}^{n-1}},
\quad t\in\R, r\in[0,+\infty[.$$
This description of the AdS metric shows  that it is conformal to a cylindrical Lorentzian  metric. In fact, we have$$
g_{AdS}=\cosh^2 r \left(-dt^2+ds^2+\sin^2 s\,\langle\;,\;\rangle_{{\mathbb S}^{n-1}}\right),$$ 
where $s=-\frac{\pi}{2}+2\arctan e^r\in [0,\frac{\pi}{2}[$, that is,$$
g_{AdS}=\cosh^2 r\left(-dt^2+g_{{\mathbb S}_+^{n}}\right),$$
where $g_{{\mathbb S}_+^{n}}$ is the unit round metric on the open 
hemisphere ${\mathbb S}_+^{n}$. From this conformal equivalence between 
Lorentzian metrics, it is straightforward to see that the conformal null infinity of the AdS spacetime 
is nothing but the cylinder $\R\times\S^{n-1}=\R\times \partial\,{\mathbb S}_+^{n-1}$ with the product
Loretzian structure $-dt^2+g_{{\mathbb S}^{n-1}}$.\\

In the second case, the level hypersurfaces of $h_a$ are umbilical
flat Euclidean spaces (horospheres) in $\H^n$. The corresponding foliation also allows
to write the hyperbolic metric $g_{\mathbb H}$ as another warped product and describe in a
corresponding manner the function $h_a$: $$
g_{\mathbb H}=dr^2+e^{2r}\,\langle\;,\;\rangle_{{\mathbb R}^{n-1}},\quad h_a=e^r,
\qquad r\in]0,+\infty[.$$
Consequently we obtain another expression for the AdS metric:$$
g_{AdS}=-e^{2r}\,dt^2+dr^2+e^{2r}\,\langle\;,\;\rangle_{{\mathbb R}^{n-1}},
\quad t\in\R, r\in]0,+\infty[.$$
In this case, the {\em height} function $h_a$ provides a different conformal equivalence for the AdS metric. Indeed, we have$$
g_{AdS}=e^{2r} \left(-dt^2+ds^2+s^2\,\langle\;,\;\rangle_{{\mathbb R}^{n-1}}\right),$$ 
where $s=e^{-r}\in [0,1[$, that is,$$
g_{AdS}=e^{2r}\left(-dt^2+g_{B^{n}}\right),$$
where  ${B}^{n} $ is the open unit disc of the Euclidean space 
${\mathbb R}^{n}$ and $g_{{B}^{n}}$ the flat metric. As in the first case, from this conformal equivalence for $g_{AdS}$, it is straightforward to see  that  the conformal null infinity of the AdS spacetime is viewed again as  $\R\times\S^{n-1}=\R\times\partial B^n$ with its standard cylindrical
Lorentzian structure $ -dt^2+g_{{\mathbb S}^{n-1}}$.\\

Then, it is clear that, for the choice of static vacuum triple $(M,g,V) = (\H^n,g_{\mathbb H},h_a)$ corresponding to the AdS spacetime, the positive function ${V}^{-1}=h_a^{-1}$ is a {\em defining function} for the hyperbolic space $\H^n$, which is the prototypical example of a {\em conformally compact} Riemannian manifold.  This simply means that the function $V^{-1}$  and the Riemannian metric $
{\overline g}= V^{-2}g=h_a^{-2}g_{\mathbb H}$
extend to an $(n+1)$-dimensional manifold with boundary ${\overline M}$ (${\mathbb S}_+^{n}$ or $B^n$) whose interior is diffeomorphic to $\H^n$ (for definitions, generalities and main facts on conformally compact Riemmanian manifolds, see    \cite{An1,An2,An3,Bi}).\\

Considering the above facts, we will say that a static vacuum $(M,g,V)$ is {\em Asymptotically Locally Anti de Sitter} (ALAdS, in short)
if its spatial slice  $(M,g)$ is a conformally compact Riemannian manifold and the positive smo\-oth function $V^{-1}$ is a defining function, that is, if there exists an $(n+1)$-dimensional manifold with boundary ${\overline M}$ such that $M$ is its
interior and  the function $V^{-1}$ extends smoothly to ${\overline M}$ in such a way that ${V^{-1}}_{|\partial M}=0$, 
$\big(dV^{-1}\big)_{|\partial M}\ne 0$, and the metric ${\overline g}=V^{-2}g$ extends to ${\overline M}$ as well (indeed, in order
to obtain the results below it would suffice to have a weaker regularity for ${\overline g}$). 
Hence, the extended metric ${\overline g}$ induces a Riemannian metric $\gamma={\overline g}_{|\partial M}$ on
the boundary at infinity $\partial M=\partial{\overline M}$ of  $(M,g)$. Note that, even if $M$ is
assumed to be connected, its boundary at infinity could be disconnected. The conformal class $[\gamma]$ of this
metric does not depend on the function $V$, but only on $(M,g)$ and this is why the $(n-1)$-dimensional conformal 
manifold $(\partial M,[\gamma])$ is
usually called the {\em conformal infinity} of $(M,g)$ (for details, see among others 
\cite{GL,L,An2,An3}). But it is clear that, in our setting, the metric $V^{-2}g_{|\partial M}$ plays a special role
within this conformal infinity $(\partial M,[\gamma])$. In fact,  the Riemannian manifold $(\partial M,
V^{-2}g_{|\partial M})$ determines the {\em conformal null infinity} ${\mathcal J}=
(\R\times\partial M,-dt^2+\gamma)$ of the static vacuum $({\mathcal V},g_{\mathcal V})=(\R\times M,-dt^2+g)$. \\

In general, a conformally compact Riemannian manifold is asymptotically negatively curved.
But if it  has also asymptotically constant scalar curvature, say $-n(n-1)$ after a suitable normalization, 
that is, if $R_g\rightarrow -n(n-1)$ when one approaches $\partial M$, then
the manifold is 
{\em Asymptotically Locally (weakly, according to \cite{Wa1}) Hyperbolic} (ALH in short). For this class of manifold we have  $K_g\rightarrow -1$  
at  infinity. As a consequence, the Ricci tensor satisfies $\Ric_g\rightarrow 
-(n-1)g$, that is, the manifold seems to be Einstein with Ricci curvature $-(n-1)$ when
one moves towards  infinity. \\

In our case, the conformally compact manifold
$(M,g)$ comes from the static vacuum $(M,g,V)$ and so,  from (\ref{scalar}), it has constant 
scalar curvature $R_g=-n(n-1)$. Thus, {\em the spatial slice of an ALAdS static vacuum is an ALH 
Riemannian manifold}. When the spatial factor $\partial M$ of the conformal null infinity 
${\mathcal J}$ of a given ALAdS static vacuum $(M,g,V)$ 
is an  $(n-1)$-dimensional hypersphere, we will say that it is {\em Asymptotically Anti de Sitter} (AAdS in short).
It is obvious that, in this case, the spatial slice $(M,g)$ is an {\em Asymptotically Hyperbolic} (AH in short) Riemannian manifold. In
this situation, P. T. Chru{\'s}ciel and M. Herzlich (for $n=3$) and X. Wang (for any $n$) generalized a former result by W. Boucher, G. Gibbons and G. Horowitz
(\cite[Theorem 4.5]{CS}, \cite[Theorem 1]{Wa2}, cfr.\! \cite{BGH}) by
proving the following uniqueness result :\\

\begin{quote} {\em Let $({\mathcal V},g_{\mathcal V})$ be
an $(n+1)$-dimensional AAdS static vacuum spacetime determined by a  triple $(M,g,V)$. Suppose
that $M$ is a spin manifold and that the conformal null infinity of
$({\mathcal V},g_{\mathcal V})$ is a unit round Lorentzian cylinder $(\R\times {\S}^{n-1},-dt^2+g_{{\mathbb S}^{n-1}})$. Then  $({\mathcal V},g_{\mathcal V})$ is the
AdS spacetime.  
}\end{quote}\

We will strengthen this result in two directions. First, we will allow the ALAdS static vacuum
triple $(M,g,V)$ to have a conformal null infinity ${\mathcal J}$ with  a non-spherical Riemannian slice  $\partial M$  and 
we will not require the metric $\gamma =V^{-2}g_{|\partial M}$ 
to have constant sectional curvature.\\

 In the spin setting,  our strategy will be to
slightly modify $V$, if necessary, in order 
to get a compactification of $(M,g)$ with non-negative scalar curvature and constant inner
mean curvature and then apply, to the corresponding induced metric on $\partial M$,  a previous estimate by X. Zhang and ourselves (see \cite{HMZ1,
HMZ2}) for the spectrum of the Dirac operator on  spin Riemannian manifolds which bound compact domains.  
Indeed, the goal of this paper is to prove the following result and 
to analyse some of its consequences (see Theorem \ref{rigidity} below).\\

\begin{quote}\label{main-th} {\em 
Let $({\mathcal V},g_{\mathcal V})$ be
an $(n+1)$-dimensional ALAdS static vacuum spacetime determined by a  triple $(M,g,V)$. Suppose
that $M$ is a spin manifold and that the spatial slice $(\partial M,V^{-2}g_{|\partial M})$ of 
its conformal null infinity ${\mathcal J}=(\R\times\partial M,-dt^2+V^{-2}g_{|\partial M})$ admits a non-trivial
 Killing or parallel spinor field. Then  $({\mathcal V},g_{\mathcal V})$ is the
AdS spacetime. In particular, there are no ALAdS static vacua whose null infinities 
have spatial slices admitting parallel spinors.  
}\end{quote}\

Note that, if $\partial M$ is spherical and  $\gamma=V^{-2}g_{|\partial M}$
is the unit round metric, then $\partial M$ admits a unique spin structure and the Riemannian
manifold $(\partial M,\gamma)$ supports a maximal number of independent Killing spinor
fields (see \cite[p.\!\! 37]{BFGK}, \cite[Examples A.1.3.2]{Gi}). Then the rigidity result
by X. Wang and P. Chru{\'s}ciel and M. Herzlich is a 
direct consequence of Theorem \ref{rigidity}.   
\\

Note that, it is well-known 
that there are compact non-spherical spin conformal manifolds carrying non-trivial
Killing and parallel spinors (see \cite{Ba1} and \cite{W1,W2}). Thus, Theorem \ref{main-th}
 implies the non-existence of ALAdS spin static vacua with these types of conformal spatial infinities.\\  

\section{Special conformal compactifications of spatial slices}
We already observed  that the spatial slices $(M,g)$ of an ALAdS static vacuum $(M,g,V)$ are ALH
manifolds. This geometrical feature reflects  an equivalent analytical property 
which is satisfied by all the defining functions $\rho$ 
of $(M,g)$, namely, $|{\overline\nabla}\rho|_{|\partial M}=1$, where the gradient and the length are
computed with respect to the same extended metric ${\overline g}=\rho^2g$ (see, for instance,
\cite[p.\! 59]{Be} or \cite[Appendix]{An3}) corresponding to the defining function $\rho$ itself. This means that the vector field ${\overline\nabla}\rho=\rho^{-2}\nabla
\rho$, resricted to $\partial M$, is an inner unit
field normal to $\partial M$ with respect to the extended metric ${\overline g}$. In particular, since $V^{-1}$ is 
a defining function for $(M,g)$, we have that$$
{{\overline\nabla}V^{-1}}_{|\partial M}=-\nabla V_{|\partial M}\perp T\partial M,\qquad |{\overline\nabla}V^{-1}|_{|\partial M}
=\left.\frac{|\nabla V|}{V}
\right|_{\partial M}=1,$$
where  the last equality is to be understood as a limit. Also, for any other defining function $\rho$ defined on ${\overline M}$,
we have
\begin{equation}\label{V-extend}
1=\left.\frac{|\nabla V|}{V}\right|_{\partial M}=\left.\frac{|\nabla\rho|}{\rho}\right|
_{\partial M},\qquad 0<(\rho V)_{|\partial M}<+\infty,
\end{equation}
where the last equality holds since the two vector fields $\nabla V$ and $\rho^{-2}{\nabla\rho}$ are parallel along $\partial M$. 
\\

In the general ALH case, where $(M,g)$ is not coming from a static vacuum, hence
not necessarily of constant scalar curvature $-n(n-1)$, if we assume 
that it is conformally compact 
of class at least $C^{3,\alpha}$, we may modify any defining function in a suitable way. 

\begin{lemma}{\cite[Section 3]{An2}, \cite[Lemma 5.4]{AD}, \cite[Lemma 5.2]{GL},  \cite[Lemma 5.1]{L}, 
\cite[Lemma 2.2]{Wa1}} \label{defining-geodesic}
Let $(M,g)$ be an ALH manifold  of class $C^{m,\alpha}$, $m\ge 3$. For each choice of a metric
$\gamma$ on its conformal infinity $\partial M$, there exists 
a unique defining function 
$r\in C^{m-1,\alpha}({\overline M})$ such that the  extended conformal metric ${\widehat g}=r^2g$ 
is of class $C^{m-1,\alpha}$, $\gamma={\widehat g}_{|\partial M}$  and $|{\widehat\nabla} r|
\equiv 1$ in a collar neighbourhood of $\partial M$.
\end{lemma}

We will say that such a function $r$ is the {\em geodesic defining function} associated with
the metric $\gamma$ in the conformal class induced on $\partial M$.  This terminology is due to the fact that, near the infinity, 
the metric $g$ of $M$ can be written as $$
g=\frac1{r^2}{\widehat g}=\frac1{r^2}(dr^2+g_r),
$$
where $g_r$ is a curve of metrics on $\partial M$ (for the regularity of $g_r$ with respect to $r$ in the Einstein case, see \cite{CDLS}). In other words, the conformal metric
${\widehat g}=dr^2+g_r$ extended to ${\overline M}$ is locally expressed  taking the ${\widehat g}$-distance to the boundary $\partial M$ 
as a first coordinate.\\

The existence of geodesic defining functions on
conformally compact ALH manifolds  is a key point to prove many of their main features.  In particular, for example the connectedness of its boundary when one
has both a lower bound on the Ricci tensor and a decay condition at infinity on the scalar curvature, along with
the  non-negativity of the conformal structure at the infinity. This connectedness was proved by E. Witten
and S.-T. Yau  for positive conformal infinities and M. Cai and G. Galloway for null conformal infinities 
(\cite{WY,CG}, \cite[Theorem 4.1]{An1}, cf.\! also \cite[Theorem 2]{HM2}). When the given ALH manifold $(M,g)$ is the spatial slice of a 
static vacuum $(M,g,V)$, on one hand, we know from (\ref{scalar}) that the scalar curvature is not
only asymptotically decaying to $-n(n-1)$, but it is in fact exactly $-n(n-1)$. However, 
on the other hand, we have no lower bounds for the Ricci tensor. So the aforementioned connectedness results do not directly apply to our situation. Yet, it is possible to make
a small detour to finally obtain a similar conclusion as in the Riemannian case.  

\begin{proposition}\label{connectedness}
Let $(M,g,V)$ be an ALAdS static vacuum triple such that the boundary at infinity $(\partial M,V^{-2}g_{|\partial M})$ of its spatial slices $(M,g)$ has  non-negative scalar curvature. Then $\partial M$ is connected.
\end{proposition}

\pf 
The conditions (\ref{vacuum}), ensuring that the Lorentzian metric $g_{\mathcal V}=-V^2dt^2+g$ on the product ${\mathcal V}=\R\times M$ is a solution of the Einstein vacuum equation $\Ric_{\mathcal V}=\Lambda\,g_{\mathcal V}$ with cosmological constant $\Lambda=-n$, also imply that the Riemannian metric ${\tilde g}=V^2d\theta^2+g$ defined on $\S^1\times M$, for a $2\pi$-periodic $\theta$, also satisfies the Einstein equation $\Ric_{\tilde g}=-n
{\tilde g}$ (see \cite{O'N}). Since $V^{-1}$ is a defining function for the conformally compact manifold
$(M,g)$, we deduce that it is also 
 a defining function for $(\S^1\times M,{\tilde g})$ with the corresponding conformal metric $V^{-2}{\tilde g}=d\theta^2+V^{-2}g$ extending to the compact manifold $\S^1\times{\overline M}$. Then $(\S^1\times M,{\tilde g})$ is a Poincar{\'e}-Einstein (PE) manifold (that is, ALH  manifolds satisfying the Einstein condition) whose conformal infinity is $(\S^1\times\partial M,[d\theta^2+V^{-2}g_{|\partial M}])$. As we are assuming that the metric $V^{-2}g_{|\partial M}$ has non-negative scalar curvature, the same occurs for the product metric $d\theta^2+V^{-2}g_{|\partial M}$ since their scalar curvatures coincide. Then the Yamabe invariant of the conformal manifold $(\S^1\times\partial M,[d\theta^2+V^{-2}g_{|\partial M}])$ is non-negative. Consequently, we can apply the aforementioned connectedness results for PE manifolds and conclude that
$\S^1\times\partial M$ must be connected. So the same conclusion is valid for $\partial M$. 

\qed\\

On the other hand, in order to study the continuous spectrum of the Laplacian on the PE manifold $(M,g)$,
 J. Lee proved in \cite{L} the existence of positive solutions of the eigenvalue equation
$\Delta u=(n+1)u$ and studied their growth at infinity $\partial M$. Some years
later, J. Qing used in \cite{Q} these particular eigenfunctions as defining functions to {\em
partially} compactify $(M,g)$ and gave
a beautiful proof, based on the positive mass theorem, of the rigidity of the hyperbolic space among all the PE manifolds 
with prescribed conformal infinity, a round sphere. In fact, he aimed to drop
the spin assumption in the corresponding result by L. Andersson and M. Dahl (\cite{AD}). \\

In our case, $(M,g)$ is the spatial slice of a static vacuum $(M,g,V)$ and thus it is an ALH manifold with
constant scalar curvature, although  not necessarily Einstein. However, we have  a
positive solution of the eigenvalue problem for the Laplace operator $\Delta$ of $g$ studied by J. Lee, namely,
the function $V$. Then we adapt 
the results  in \cite[Section 5]{L} in order to get, by  Lemma \ref{defining-geodesic}, 
a sufficiently good description
of the asymptotical behaviour of $V$ in terms of the geodesic defining function on $(M,g)$ associated with
$V^{-1}$. For this, we need to prove the following two  lemmae. 

\begin{lemma}\label{R^*}
Let $(M,g,V)$ be an ALAdS static vacuum and $\varepsilon\ge 0$ a real number. Then the function $(V+\varepsilon)^{-1}$
is a defining function for the ALH manifold $(M,g)$ such that the corresponding extended Riemannian metric $g^*=(V+\varepsilon)^{-2}g$
has scalar curvature
\begin{equation}\label{estimate*}
R_{g^*}=n(n-1)\big(V^2-|\nabla V|^2-\varepsilon^2\big).
\end{equation}
In particular, the function $V^2-|\nabla V|^2$ extends to the compact manifold with boundary ${\overline M}$.
\end{lemma}

\pf
Since $V^{-1}$ is a positive smooth function on $M$ extendable to ${\overline M}$ and vanishing on $\partial M$, we have that 
$(V+\varepsilon)^{-1}=V^{-1}/(1+\varepsilon V^{-1})$ is another positive smooth function on $M$
extendable to ${\overline M}$ and vanishing on $\partial M$ as well.  It is immediate to see that
 $$
 \big(d(V+\varepsilon)^{-1}\big)_{|\partial M}=(dV^{-1})_{|\partial M}\ne 0.$$ 
 Moreover the metric
 \begin{equation}\label{g^*}
g^*=(V+\varepsilon)^{-2}g=\left(\frac1{1+\varepsilon V^{-1}}\right)^2V^{-2}g
\end{equation}
defined on $M$  extends to a metric on ${\overline M}$ with the same regularity as that of the metric
$V^{-2}g$.  So, $(V+\varepsilon)^{-1}$ is a defining function for $(M,g)$. Moreover, 
this metric $g^*$, restricted to $\partial M$, gives  $$
{g^*}_{|\partial M}=\big((V+\varepsilon)^{-2} g\big)_{|\partial M}=(V^{-2} g)_{|\partial M}.$$
In particular, $\big({M}^*,g^*\big)=\big({\overline M},(V+\varepsilon)^{-2}g\big)$ is a
conformal compactification of the ALH manifold $(M,g)$.
Now, rewrite the conformal change between the metrics $g$ and $g^*$ as $$
g=(V+\varepsilon)^2g^*,$$
then the relation between the associated Ricci tensors $\Ric_{g^*}$ and $\Ric_g$
on the open manifold $M$ (see  \cite[p.\! 59]
{Be} or \cite[Appendix]{An3}),  is given by $$ 
 \Ric_{g^*}=\Ric_g+(n-2)\frac{\nabla^2 V}{V+\varepsilon} + 
 \frac{\Delta V}{V+\varepsilon}g-(n-1)\frac{|\nabla V|^2}{(V+\varepsilon)^2}g.
$$
Taking traces with respect to $g$ and multiplying by $(V+\varepsilon)^2$, one gets for 
the corresponding scalar curvatures$$
R_{g^*}=(V+\varepsilon)^2R_g+2(n-1)(V+\varepsilon)\Delta V-n(n-1)|\nabla V|^2.
$$
As we pointed out in (\ref{vacuum}) and (\ref{scalar}), since $R_g=-n(n-1)$ and
the function $V$ is an eigenfunction of $\Delta$ associated with the
eigenvalue $-n$, we finally
get (\ref{estimate*}).

\qed

An important consequence of  Lemma \ref{R^*} is the following information about
the asymptotical behaviour of the Ricci tensor of the spatial slice $(M,g)$ when
one approaches to $\partial M$ along the direction $\nabla V$.\

\begin{lemma}\label{Ricci}
Let $\rho$ be any defining function of the ALH spatial slice $(M,g)$ of an ALAdS static
vacuum. Then
\begin{equation}\label{Ri-inf}
\lim_{\rho\rightarrow 0}\frac1{\rho^2}\big(\hbox{\rm Ric}_g+(n-1)g\big)\left(
\frac{\nabla\rho}{\rho},\frac{\nabla\rho}{\rho}\right)=0.
\end{equation}
\end{lemma} 

\pf
From Lemma \ref{R^*}, we know that $V^2-|\nabla V|^2$ is a differentiable function on
the compact manifold ${\overline M}$. Let ${\overline N}$ be the inner unit vector field normal 
to $\partial M$ with respect to the extended metric ${\overline g}=V^{-2}g$. Then, the Lee derivative of the function $V^2-|\nabla V|^2$ in the direction
of the vector field $\frac1{2}{\overline N}$, is given by  $$
\frac1{2}{\overline N}\cdot\big( V^2-|\nabla V|^2\big)=\left(Vg(\nabla V,{\overline N})-(\nabla^2 V)
({\overline N},\nabla V)\right)_{|\partial M}.
$$
But, by (\ref{V-extend}) and previous comments, we have ${\overline N}=-\nabla V_{|\partial M}$. Thus$$
\frac1{2}{\overline N}\cdot\big( V^2-|\nabla V|^2\big)=\left((\nabla^2 V)
(\nabla V,\nabla V)-V|\nabla V|^2\right)_{|\partial M}.
$$
Now, from the first equality in (\ref{vacuum}), it follows $$
\frac1{2}{\overline N}\cdot\big( V^2-|\nabla V|^2\big)=\left(V\big(\hbox{\rm Ric}_g+(n-1)g\big)
(\nabla V,\nabla V)\right)_{|\partial M}.
$$
Using again (\ref{V-extend}), we know that$$
\nabla V_{|\partial M}=\left.\frac{\nabla \rho}{\rho^2}\right|_{\partial M},\qquad
V{\rho}=O(1),$$
hence the existence of the following limit$$
 \lim_{\rho\rightarrow 0}\frac1{\rho^3}\big(\hbox{\rm Ric}_g+(n-1)g\big)\left(
\frac{\nabla\rho}{\rho},\frac{\nabla\rho}{\rho}\right),$$
which implies (\ref{Ri-inf}).


\qed

Now, we dispose of all the necessary tools to control the asymptotic behavior of the
function $V$ in the given static vacuum $(M,g,V)$. 

\begin{proposition}{\rm (Cfr.\! \cite[Proposition 4.1, Lemmae 5.1 and 5.2]{L} and \cite[Lemma 2.1]{GQ})} \label{eigenfunctions}
Let $(M,g,V)$ be an $(n+1)$-dimensional, $n\ge 3$, ALAdS static vacuum  and  $r$ the geodesic defining function associated with 
the metric $\gamma=V^{-2}g_{|\partial M}$ on the conformal   infinity $\partial M$ of its spatial slices, according to Lemma \ref{defining-geodesic}. 
Then, one has
\begin{equation}\label{u-r}
V=\frac1{r}+hr,
\end{equation}
for a function
 $h\in C^{2,\alpha}(M)\cap C^0({\overline M})$ such that
\begin{equation}\label{asymptotic-h}
h_{|\partial M}=\frac{R_\gamma}{4(n-1)(n-2)},\qquad |\nabla h|=O\big(r^{\frac{\alpha}{2}}\big),
\end{equation}
where $\alpha>0$ and $R_\gamma $ is the scalar curvature of the  metric $\gamma$ on 
$\partial M$.
\end{proposition}

\pf
This result relies on  some of the assertions in Proposition 4.1, Lemma 5.1 (see Lemma
\ref{defining-geodesic} of this paper) and Lemma 5.2 in \cite{L}. These assertions, 
only require $(M,g)$ to be ALH, a condition satisfied by 
$(M,g)$ since it is the spatial slice of an ALAdS static vacuum. The proof of the second of these lemmae in \cite{L} 
definitely uses the hypothesis that $(M,g)$ is Einstein. So in fact, we show now that the result holds under weaker hypotheses. \\

A careful look  at the proof of  Lemma 5.2 in \cite{L}, yields to the observation that the Einstein condition on the metric $g$ is used to obtain the equalities labelled (5.3) and (5.7) in \cite{L}, namely$$
\left\{\begin{array}{ll}
\hbox{\rm (5.3)}\quad{\widehat\Delta}r=-\frac1{2(n-1)}rR_{\widehat g}\quad\hbox{on a neighbourhood of $\partial M$},
\\ \\ \hbox{\rm (5.7)}\quad R_\gamma=\frac{n-2}{n-1}{R_{\widehat g}}_{|\partial M},
\end{array}\right.
$$
where ${\widehat g}=r^2g$ is the conformal metric extended to ${\overline M}$, $\gamma=
{\widehat g}_{|\partial M}$  its restriction to the conformal infinity and ${\widehat\Delta}$ its scalar Laplacian.
To get 
 (5.3),  consider the relation between the Ricci tensors $\Ric_{\widehat g}$ and $\Ric_g$
of the conformal metrics $\widehat g$ and $g$ on the open manifold $M$ (see again \cite[p.\! 59]
{Be} or \cite[Appendix]{An3}): 
$$
 \Ric_{g}=\Ric_{\widehat g}+(n-2)\frac{{\widehat \nabla}^2 r}{r} + 
 \frac{{\widehat\Delta} r}{r}{\widehat g}-(n-1)\frac{|{\widehat\nabla} r|^2}{r^2}{\widehat g},
$$
where ${\widehat\nabla}$ and ${\widehat\nabla}^2$ are respectively the gradient and the Hessian operator of
${\widehat g}$.
Multiplying by $r$ and putting $|{\widehat\nabla}r|^2=1$, we obtain
\begin{equation}\label{step1-r}
r\big(\hbox{\rm Ric}_g+(n-1)g\big)=r\hbox{\rm Ric}_{\widehat g}+(n-2){\widehat\nabla}^2r+
({\widehat\Delta}r){\widehat g} 
.
\end{equation} 
Taking traces with respect to ${\widehat g}$, knowing that  
 $R_g=-n(n-1)$,  it follows that on a collar neighbourhood of $\partial M$ one has
\begin{equation}\label{step2-r}
rR_{\widehat g}+2(n-1){\widehat\Delta}r=\frac1{r}\big(R_g+n(n-1)\big)=0,
\end{equation}
giving (5.3) in \cite{L}. \\

As for the second, 
  putting first $r=0$ in (\ref{step2-r}), then dividing this same expression by $r$ and
taking limits as $r\rightarrow 0$, we get
\begin{equation}\label{step4-r}
({\widehat\Delta r})_{|\partial M}=0,\qquad ({\widehat\nabla}r\cdot{\widehat\Delta}r)_{|\partial M}
=-\frac1{2(n-1)}{R_{\widehat g}}_{|\partial M}.
\end{equation}
On the other hand, in the collar neighbourhood of $\partial M$, taking derivatives of $|{\widehat \nabla}r|^2= 1$
with respect to a
vector field $X$ tangent to ${\overline M}$, yield  $({\widehat\nabla}^2r)({\widehat\nabla}r,X)=0$. In particular,
\begin{equation}\label{hessian-r}
({\widehat\nabla}^2 r)({\widehat\nabla}r,{\widehat\nabla}r)=0.
\end{equation}
Again, taking derivatives  of $({\widehat\nabla} dr)({\widehat\nabla}r,X)=0$ with respect to the same direction $X$,
using the Ricci equation and then taking traces, imply
\begin{equation}\label{laplacian-r}
|{\widehat\nabla}dr|^2+{\widehat\nabla}r\cdot{\widehat\Delta}r
+\hbox{\rm Ric}_{\widehat g}({\widehat\nabla}r,{\widehat\nabla}r)=0.
\end{equation} 
 Now, we apply the tensorial equality (\ref{step1-r}) 
to the vector field ${\widehat\nabla}r$ and use (\ref{hessian-r}), to get
\begin{equation}\label{step3-r}
\frac1{r}\big(\hbox{\rm Ric}_g+(n-1)g\big)\left(\frac{\nabla r}{r},\frac{\nabla r}{r}\right)=r
\hbox{\rm Ric}_{\widehat g}\big({\widehat\nabla}r,{\widehat\nabla}r\big)+
{\widehat\Delta}r.
\end{equation}
Dividing also (\ref{step3-r}) by $r$, taking limits as $r\rightarrow 0$, keeping in mind the
first equality in (\ref{step4-r}) and using Lemma \ref{Ricci}, we 
obtain $$
({\widehat\nabla}r\cdot{\widehat\Delta}r)_{|\partial M}=
-\hbox{\rm Ric}_{\widehat g}\big({\widehat\nabla}r,{\widehat\nabla}r\big)_{|\partial M}.$$ 
This together with (\ref{laplacian-r}) and the second equality in (\ref{step4-r}) give $$
{\widehat \nabla}^2r_{|\partial M}=0,\qquad \hbox{\rm Ric}_{\widehat g}({\widehat
\nabla}r,{\widehat\nabla}r)_{|\partial M}=\frac1{2(n-1)}{R_{\widehat g}}_{|\partial M}.$$\

Now, take into account that the Hessian ${\widehat\nabla}^2r$, restricted to directions
orthogonal to the gradient ${\widehat\nabla}r$, is the opposite of the second fundamental form
${\widehat\sigma}$ of the level hypersurfaces $r=r_0$ with respect to the choice of inner unit normal  
${\widehat N}={\widehat\nabla}r$ and to the metric ${\widehat g}$ on ${\overline M}$.
So, using (\ref{hessian-r}), we see that ${\widehat H}=-\frac1{n-1}{\widehat\Delta}r$ is the
mean curvature function of these level hypersurfaces. Using such submanifold theory language, we may
rewrite the last relations as 
\begin{equation}\label{H^}
{\widehat\sigma}_{|\partial M}=0,\quad {\widehat H}_{|\partial M}=0,\quad
\hbox{\rm Ric}_{\widehat g}({\widehat
N},{\widehat N})_{|\partial M}=\frac {{R_{\widehat g}}_{|\partial M}}{2(n-1)}.
\end{equation}
Knowing that the Gau{\ss} equation relating the scalar curvature $R_\gamma$ of the metric $\gamma$
induced on the boundary $\partial M$ and the restriction ${R_{\widehat g}}_{|\partial M}$ is given by
$$
R_\gamma={R_{\widehat g}}_{|\partial M}-2\hbox{\rm Ric}_{\widehat g}({\widehat N},
{\widehat N})_{|\partial M}+(n-1)^2{\widehat H}^2_{|\partial M}-|{\widehat\sigma}
|^2_{|\partial M},$$
the relations in  (\ref{H^}) imply  equality  (5.7) in \cite{L}.

\qed

As mentioned before, the idea of using the eigenfunctions of the Laplacian, with controlled behaviour at infinity, as
defining functions on a given PE space $(M,g)$ is due to J. Qing 
(see \cite{Q}), although he conformally modifies the original complete manifold through these eigenfunctions 
without actually compactifying it. Instead, he gets a {\em partial} compactification, that is, a conformal complete Riemannian manifold $({\overline M},{\overline g})$ 
whose boundary $\partial M$  is diffeomorphic to the Euclidean space ${\mathbb R}^n$ and 
such that the Riemannian manifold constructed by doubling $({\overline M},{\overline g})$ along this 
boundary is an asymptotically Euclidean manifold without boundary and with non-negative integrable scalar curvature. Then a suitable 
use of the positive mass theorem allows him to go on with his reasoning. \\

Since our ALH manifold $(M,g)$ is the spatial slice of a static vacuum $(M,g,V)$, we know from (\ref{scalar}) that $V$ itself is such an eigenfunction of the Laplacian. Moreover, Proposition 
\ref{eigenfunctions} gives a reasonable control on $V$ near the conformal spatial infinity of
$(M,g)$. As in Lemma \ref{R^*}, we will proceed  by slightly modifying $V$  in order to {\em totally} compactify $(M,g)$.

\begin{theorem}\label{compactification}
Let $(M,g,V)$ be an $(n+1)$-dimensional, $n\ge 3$, ALAdS static vacuum and suppose that the spatial
slice  $(\partial M,\gamma=V^{-2}g_{|\partial M})$ of its conformal null infinity has non-negative  scalar curvature.
 Then there exists a  defining 
function $\rho^*$ for $(M,g)$ such that, 
if $g^*$ is the extension of $(\rho^*)^2g$ to ${\overline M}$, $$
{g^*}_{|\partial M}=\gamma,\qquad R_{g^*}\ge 0,\qquad H^*=\varepsilon,$$
where  $R_{g^*}$ is the scalar curvature of the compact Riemannian manifold $({\overline
M},g^*)$,  
$H^*$ is the (inner) mean curvature  of the conformal infinity $\partial M$ as a hypersurface
of $({\overline M},g^*)$ and $\varepsilon\ge 0$ is given by$$
\varepsilon=\inf_{\partial M}\sqrt{\frac{R_\gamma}{(n-1)(n-2)}}. $$ 
\end{theorem}

\pf 
 As in Lemma \ref{R^*},
we define the function $\rho^*=(V+\varepsilon)^{-1}$ for this precise choice of  $\varepsilon$. In  Lemma \ref{R^*} we also showed that $\rho^*$ is indeed a defining function for $(M,g)$ with ${g^*}_{|\partial M}=\gamma$ and established that   the scalar curvature of the compact Riemannian with boundary $({\overline M},g^*)$, satisfies (\ref{estimate*}).  

Now, to show that $R_{g^*}\ge 0$, our approach  is to study the points where the function $\Phi :=V^2-|\nabla V|^2-\varepsilon^2$, defined on
${\overline M}$,  attains its minimum. Taking into account the Bochner formula for the Laplacian of
the squared length of a gradient, we have$$
\frac1{2}\Delta \Phi=|\nabla V|^2+V\Delta V-|\nabla^2V|^2-
g(\nabla V,\nabla\Delta V)-\Ric_g(\nabla V,\nabla V).$$ 
Using (\ref{vacuum}), we obtain
\begin{equation}\label{Bochner}
\frac1{2}\Delta \Phi= -|\nabla^2V-Vg|^2-\frac{(\nabla^2V)(\nabla V,\nabla V)}{V}+|\nabla V|^2.
\end{equation}
On the other hand, from the  definition of $\Phi$, we have
\begin{equation}\label{gradient-Fi}
g(\nabla\Phi,\nabla V)=V|\nabla V|^2-(\nabla^2 V)(\nabla V,\nabla V).
\end{equation}
If the minimum of $\Phi$ is attained at a point $p\in M$, then (\ref{gradient-Fi}) and (\ref{Bochner}) would imply that $(\Delta\Phi)(p)\le 0$ and thus the strong minimum principle implies that $\Phi$ is constant. Consequently, the minimum of $\Phi$ could be also attained
at the infinity $\partial M$. Hence, we can assume that $\Phi$ reaches its minimum
value at $\partial M$.    
\\

To study the asymptotical behaviour of the function $\Phi$, we can use the expression (\ref{u-r}) of the eigenfunction $V$ in terms of the 
geodesic defining function $r$. Then$$
\nabla V=\left(h-\frac1{r^2}\right)\nabla r+r\nabla h.$$
Taking squared norms with respect to the metric $g$, we get$$
|\nabla V|^2=\left(h-\frac1{r^2}\right)^2|\nabla r|^2+r^2|\nabla h|^2+2r\left(h-\frac1
{r^2}\right) g(\nabla r,\nabla h).$$ 
But we know that $|\nabla r|^2=r^2$ near $\partial M$ due to the geodesic 
character of the defining function $r$. 
Putting this into the last equation and using (\ref{u-r})
again, we have$$
\Phi=4h-r^2|\nabla h|^2-2r\left(h-\frac1{r^2}\right)g(\nabla r,\nabla h),$$
which is valid in a collar neighbourhood of $\partial M$.
From (\ref{asymptotic-h}), we know that $h$ extends to $C^0({\overline M})$ 
and  that $|\nabla h|=O\big(r^\frac{\alpha}{2}\big)$. So, the third term 
in the right side of the previous equation  satisfies$$
\left|r\left(h-\frac1{r^2}\right)g(\nabla r,\nabla h)\right|\le \left|hr^2-1\right||\nabla h|=
O\big(r^\frac{\alpha}{2}\big)$$
as a consequence of the Schwarz inequality for $g$ and again from the fact that $|\nabla r|=r$ in a collar
neighbourhood of $\partial M$.  From this inequality and the equality above, we conclude, taking limits when $r\rightarrow 0$, that$$
\Phi_{|\partial M}=\left(V^2-|\nabla V|^2\right)_{|\partial M}=4h_{|\partial M}
=\frac{R_\gamma}{(n-1)(n-2)}\ge\varepsilon^2,$$
where we have used again (\ref{asymptotic-h}) and the choice of $\varepsilon$. \\

Thus $\Phi_{|\partial M}\ge\varepsilon^2$ together with equality (\ref{estimate*})
imply that on $ {\overline M}$,   one has 
\begin{equation}\label{positivity-u}
R_{g^*}=n(n-1)\big(V^2-|\nabla V|^2-\varepsilon^2\big)\ge 0 .
\end{equation}

To finish the proof, it remains to compute the mean curvature $H^*$ of the conformal
infinity $\partial M$ as a hypersurface of the compactified Riemannian manifold
$({\overline M},g^*)$. Observe that, by definition of $\rho^*$ and Proposition \ref{eigenfunctions},$$
g^*=(\rho^*)^2g=\left(\frac{r}{1+\varepsilon r+hr^2}\right)^2g=\left(\frac1{1+\varepsilon r+hr^2}\right)^2
{\widehat g}$$
where ${\widehat g}=r^2g$ is the extended metric on ${\overline M}$ corresponding to the
geodesic defining function $r$. 
But,  we know from (\ref{H^}) that  the mean curvature ${\widehat H}$ of the spatial conformal  infinity 
${\partial M}$ as a hypersurface of $({\overline M},{\widehat g})$ vanishes (in fact, it is
 a totally geodesic hypersurface). 
So, in order to compute $H^*$, it suffices to use the well-known relation between the two mean curvatures of a hypersurface corresponding to two conformal metrics on the ambient space  (see, for instance, \cite{E} or \cite[(4.4)]{HMZ2}):$$
H^*=\frac1{R}\big({\widehat H}-{\widehat g}( {\widehat\nabla}\log R,{\widehat N})\big)=-\frac1{R^2}
{\widehat g}({\widehat\nabla}R,{\widehat N}),$$
where $R=(1+\varepsilon r+
hr^2)^{-1}$ and ${\widehat N}$ is the inner unit normal  along $\partial M$ with respect to the
metric ${\widehat g}$, which can be also written as 
${\widehat\nabla} r$. Since $$
{\widehat\nabla} R_{|\partial M}=-\left.\frac{\varepsilon{\widehat\nabla} r+
{\widehat\nabla} (r^2h)}
{\big(1+\varepsilon r+hr^2\big)^2}\right|_{\partial M}=-\varepsilon{\widehat\nabla} r_{|\partial M},$$  
we finally obtain$$
(H^*)_{|\partial M}=
\varepsilon{\widehat g}({\widehat\nabla} r, 
{\widehat\nabla} r)_{|\partial M}=\varepsilon.$$
\qed

As a first consequence, we have that our assumption on the scalar curvature of
 $(\partial M,\gamma)$  controls the relative Yamabe invariant  (see \cite{E} for a definition)
of the  compactifications of the spatial slice $(M,g)$ of the static vacuum $(M,g,V)$.

\begin{corollary}
Let $(M,g,V)$ be an $(n+1)$-dimensional ALAdS static vacuum with $n\ge 3$ and whose  conformal spatial infinity 
$(\partial M,V^{-2}g_{|\partial M})$ has 
positive (respectively non-negative) scalar curvature. Then the relative
Yamabe invariant of the conformal compactification $({\overline M},{\overline g}=V^{-2}g)$
is also positive (respectively non-negative).
\end{corollary}

\pf 
In the context of the Yamabe problem for $n$-dimensional compact manifolds with boundary, J. Escobar introduced
in \cite{E} the following eigenvalue problem$$
\left\{\begin{array}{ll}
-\frac{4(n-1)}{n-2}{\overline\Delta}f+R_{\overline g}f=0,\quad\hbox{on }{\overline M},\\ \\
-\frac{2}{n-2}{\overline g}({\overline\nabla}f,{\overline N})+{\overline H}f={\overline \nu} f,\quad\hbox{along }\partial M,
\end{array}\right.
$$
for functions $f\in C^1({\overline M})$. He proved that the sign of the first eigenvalue
${\overline\nu}_1$ of this problem (if finite) is a conformal invariant of the metric ${\overline g}$
and that its sign coincides with that of the relative Yamabe invariant of the
manifold with boundary ${\overline M}$, whose value is a conformal invariant as well
(see \cite{HMZ2} for a relation between $\nu_1$, the relative Yamabe invariant and the Dirac operator of the boundary).
Under our hypotheses, we can apply Theorem \ref{compactification} 
and dispose of a metric $g^*$ on ${\overline M}$ conformal to ${\overline g}$ and such that
$R_{g^*}\ge 0$ on ${\overline M}$ and $H^*=\varepsilon$ along $\partial M$, where $\varepsilon
>0$ (respectively $\varepsilon\ge0$) according to whether the scalar curvature of $\partial M$ is positive (respectively
non-negative). Since the metrics ${\overline g}$ and $g^*$ are conformal we can compute $\nu_1^*$ in order
to determine the sign of ${\overline\nu}_1$. For this, take a non-trivial function $f\in C^1({\overline
M})$. We have$$
\int_{\overline M}\left(\frac{2}{n-1}|\nabla^* f|^2+\frac1{2n}R_{g^*}f^2\right)+
\int_{\partial M}H^*f^2\ge \varepsilon\int_{\partial M}f^2.$$
Hence, by the variational characterization of the eigenvalue $\nu_1^*$, we have $\nu_1^*
\ge \varepsilon$ and so $\nu_1^*$ (and hence ${\overline \nu}_1$) is positive when $\varepsilon>0$
and non-negative when $\varepsilon\ge0$.

\qed

\section{Spatial slices admitting Killing or parallel spinors on its infinity}
Suppose now that the spatial slice $M$ of an ALAdS static vacuum $(M,g,V)$ is a spin manifold on which we fix 
 a spin structure. We know that its spatial slice
$(M,g)$ is an ALH  manifold. Indeed, the metric ${\overline g}=V^{-2}g$ on $M$ extends to 
a compact manifold with boundary ${\overline M}$ whose interior is $M$ itself. Then, it is 
straightforward to check that ${\overline M}$ is also a spin manifold and that we may fix a unique spin structure
on ${\overline M}$ such that its restriction on the open subset $M$ is precisely the given spin
structure of $M$. Since the spatial conformal 
infinity $\partial M$ is always an orientable hypersurface (recall that the gradient of a geodesic
defining function provides a global unit normal field), we have that ${\partial M}$ is also a 
spin manifold and   that an induced spin structure on the infinity is inherited from the fixed structure  on
${\overline M}$. Moreover, for the Riemanian metric ${\overline g}$ on ${\overline M}$ we have 
an associated spinor bundle $({\mathbb S}{\overline M},{\overline\nabla },{\overline {\mathfrak c}})$,
where ${\overline\nabla}$ is the spin Levi-Civita connection and ${\overline{\mathfrak c}}$ is the Clifford
multiplication (for generalities on spin structures see any of \cite{BFGK,BHMM,Fr,Gi,LM}). It is a well-known
fact that the restriction to the hypersurface $\partial M$ of the spinor bundle ${\mathbb S}{\overline M}$ 
can be identified with one or two copies of the spinor bundle corresponding to the induced spin structure and 
the induced Riemannian metric $\gamma={\overline g}_{|\partial M}$ according to the parity of the dimension $n$
of ${\overline M}$. More precisely, if $\varphi$ is a section
of the restricted bundle ${\mathbb S}{\overline M}_{|{\partial M}}$, we consider the new Clifford multiplication$$
{\mathfrak c}^{\partial M}(X)\varphi={\overline{\mathfrak c}}(X){\overline{\mathfrak c}(N)}
\varphi.
$$
and the new connection$$
\nabla^{\partial M}_X\varphi={\overline\nabla}_X\varphi-\frac1{2}{\overline{\mathfrak c}}({\overline A}X){\overline{\mathfrak c}
({\overline N})}
\varphi=
{\overline\nabla}_X\varphi-\frac1{2}{\mathfrak c}^{\partial M}({\overline A}X)\varphi,\quad\forall X\in\Gamma(T\partial M),$$
where ${\overline N}$ is the (inner) unit normal field along $\partial M$ and ${\overline A}$ is its corresponding shape
operator. Then, we have an isomorphism$$
({\mathbb S}{\overline M}_{|\partial M},\nabla^{\partial M},{\mathfrak c}^{\partial M})\cong
\left\{\begin{array}{ll}
({\mathbb S}\partial M,\nabla,{\mathfrak c}),\hbox{ if $n$ is odd} \\ \\
({\mathbb S}\partial M,\nabla,{\mathfrak c})\oplus({\mathbb S}\partial M,
\nabla,-{\mathfrak c}),
\hbox{ if $n$ is even},
\end{array}\right.$$  
 where $({\mathbb S}\partial M,\nabla,{\mathfrak c})$ is the spinor bundle corresponding to the 
 spin structure and to the Riemannian metric induced on $\partial M$ (for this relationship between the spinor bundles on
 a hypersurface and on its ambient space, see, for instance, \cite{Ba3,Fr, HM1,HMZ1,HMZ2}). 
 Due to this identification we can say that each spinor field on 
 ${\overline M}$ determines, by restriction, a spinor field on the boundary $\partial M$
 and we can talk about possible extensions to ${\overline M}$ of the spinor fields defined
 on $\partial M$.\\
 
 Let $\varphi\in \Gamma({\mathbb S}\partial M)$ be a spinor field on the Riemannian slice
 of the conformal
 null infinity of an ALAdS spin static vacuum $(M,g,V)$.  When $\varphi$  satisfies the first order  equation
 \begin{equation}\label{Killing-spinor}
 \nabla_X\varphi+\frac{\lambda}{n-1}{\mathfrak c}(X)\varphi=0,\quad \lambda\in{\mathbb C},\quad\forall X\in
 \Gamma(T\partial M),
 \end{equation}
 we will say that $\varphi$ is a {\em Killing spinor} if $\lambda\in{\mathbb C}^*$. Of course, if the same equation
 is satisfied with $\lambda=0$, we say that $\varphi$ is a {\em parallel spinor} (we refer  to \cite{BFGK,
 BHMM,CGLS,Gi} for definitions and main properties). It can be shown that $\lambda$ has to be
  a real or a purely imaginary number. So, we will talk about {\em  real or imaginary 
 Killing spinors} according to $\lambda\in{\mathbb R}^*$ or $\lambda\in i{\mathbb R}^*$.
 Of course, these
 definitions are usually made for general spin Riemannian  manifolds, not necessarily boundaries
 of compactifications. It is immediate that a Killing or parallel spinor must be an 
 eigenspinor for  the well-known Dirac operator
$D$ locally defined by$$
 D\varphi=\sum_{i=1}^{n-1}{\mathfrak c}(e_i)\nabla_{e_i}\varphi,$$
 where $\{e_1,\dots,e_{n-1}\}$ is a local orthonormal frame on $\partial M$. In fact, 
  (\ref{Killing-spinor}) immediately implies $D\varphi=\lambda\varphi$.
 \\
 
 The existence of parallel or Killing spinors imposes strong restrictions on the geometry of the manifold
 and on its holonomy. Such manifolds have to be Einstein with scalar curvature $R=4(n-1)(n-2)\lambda^2$.
 Indeed, M. Wang, H. Baum and C. B{\"a}r (\cite{W1,W2,Bm,Ba1,Ba2})  classified some types of spin Riemannian
 manifolds admitting non-trivial parallel, imaginary Killing and  real Killing spinors,
 respectively.  When the considered spin Riemannian manifold is compact,
 since the eigenvalues of its Dirac operator have to be real, Killing spinors must be real as well and,
 moreover, as it was shown by T. Friedrich \cite[Corollary 1, Theorem 9]{BFGK}, they are eigenspinors corresponding to the eigenvalues 
 with the least absolute value $\pm\sqrt{\frac{(n-1)R}{4(n-2)}}$. This quick review about spin structures and Killing spinors allows us to set up our first rigidity result.\\

 \begin{theorem}\label{rigidity}
  Let $({\mathcal V},g_{\mathcal V})$ be
an $(n+1)$-dimensional, $n\ge 3$, ALAdS static vacuum spacetime determined by a static vacuum triple $(M,g,V)$. Suppose
that $M$ is a spin manifold and that the spatial slice $(\partial M,\gamma)$ of its conformal null infinity ${\mathcal J}=(\R\times\partial M,-dt^2+\gamma)$, with $\gamma= V^{-2}g_{|\partial M}$, admits a non-trivial
Killing or parallel spinor field. Then  $({\mathcal V},g_{\mathcal V})$ is the
AdS spacetime. Consequently, the parallel case cannot occur.
 \end{theorem} 
 
\pf 
Let $(M,g,V)$ be a static vacuum triple determining $({\mathcal V},g_{\mathcal V})$ and let $({\overline M},{\overline g}
=V^{-2}g)$ be the conformal spin compactification of its spatial slice $(M,g)$. Take a non-trivial Killing or parallel spinor $\varphi\in\Gamma({\mathbb S}\partial M)$ on the spatial
slice $(\partial M,\gamma)$ of the conformal
null infinity ${\mathcal J}$. Since $\partial M$ is compact, we know that $\lambda$ in  (\ref{Killing-spinor}) is a real number
and that the metric $\gamma$ on $\partial M$ has 
 constant scalar curvature $R_{\gamma}=(n-1)(n-2)\lambda^2$. So, we have on $\partial M$ either positive constant scalar curvature and a non-trivial 
 real Killing spinor $\varphi_1$,
or identically zero scalar curvature and a non-trivial parallel spinor $\varphi_0$. In particular,  we can apply Theorem
\ref{compactification} to this situation by choosing $\varepsilon=|\lambda|$. So we have a metric $g^*$ on the spin compactification 
${\overline M}$  such that $R_{g^*}\ge 0$, $H^*=\varepsilon$
and ${g^*}_{|\partial M}=\gamma$. Under these conditions, we may use a lower estimate for 
the spectrum of the Dirac operator $D$ of $(\partial M,\gamma)$ obtained by X. Zhang and the authors in \cite{HMZ1} (see \cite[Theorem 3.7.1]
{Gi}). This result asserts that, if $\lambda_1(D)$ stands for the eigenvalue of $D$ with the lowest absolute value, 
then $$
|\lambda_1(D)|\ge \frac{n-1}{2}\varepsilon$$ 
and, if the equality holds, then the eigenspace
associated with $\lambda_1(D)$ is built from  parallel spinor fields 
on $({\overline M}, g^*)$. But we know that there exists on $\partial M$
a non-trivial Killing or parallel spinor $\varphi_{i(\varepsilon)}$, with $i(\varepsilon)=1$ or 
$i(\varepsilon)=0$, depending on whether
$\varepsilon>0$ or $\varepsilon=0$ and hence 
$D\varphi_{i(\varepsilon)}=\pm
\frac{n-1}{2}\varepsilon\varphi_{i(\varepsilon)}$. Then the equality $|\lambda_1(D)|=\frac{n-1}{2}
\varepsilon$ holds and so $\varphi_{i(\varepsilon)}$ comes from  a  parallel spinor field
$\Psi\in{\mathbb S}{\overline M}$. Note that
$\Psi$ has to be a non-trivial parallel spinor since its restriction to $\partial M$ is non-trivial.\\

It was shown by Hitchin in \cite{Hit} (see also \cite[Chapter 6]{BFGK}) that the existence
of a non-trivial parallel spinor forces the Ricci tensor to vanish identically. Then $\Ric_{g^*}
=0$ on ${\overline M}$ and so $R_{g^*}=0$ as well. From (\ref{estimate*}), 
(\ref{Bochner}), (\ref{gradient-Fi}) and (\ref{positivity-u}), and the proof of Theorem \ref{compactification}, we conclude 
\begin{equation}\label{Obata}
\nabla^2V=Vg,\qquad V^2-|\nabla V|^2-\varepsilon^2=0.
\end{equation}
Hence the complete manifold $(M,g)$ admits a
non-trivial (in fact, positive) solution $V$ to the Obata type equation $\nabla^2V=Vg$.
If the function $V$ has a critical point, Theorem C in \cite{K}  implies that $(M,g)$ is
isometric to the hyperbolic space $(\H^n,g_{\mathbb H})$ and $V$ is a positive {\em height} function.
Then $({\mathcal V},g_{\mathcal V})$ is isometric to the  AdS spacetime.
Assume on the contrary that $V$ has no critical points on $M$. We can normalize
the gradient $\nabla V$ to obtain a global unit vector field $X=\frac{\nabla V}
{|\nabla V|}$ on $M$ satisfying $\nabla_XX=0$ as a consequence of (\ref{Obata}). Hence
the integral curves of $X$ are geodesics and are defined on the whole real line. Take
a positive real number $a$ in the image of $V$. Then $P=V^{-1}(\{a\})$ is a closed hypersurface 
in $M$ and so compact because $V$ tends to $+\infty$ when one approaches $\partial M$
(see (\ref{u-r})). 
Let ${\mathcal F}:M\times {\mathbb R}
\rightarrow M$ be the flow of $X$. From the considerations above it follows that the
restriction ${\mathcal F}: P\times {\mathbb R}\rightarrow M$ is a diffeomorphism with$$
{\mathcal F}(p,s)=\gamma_p(s),\qquad\forall p\in P,\, \forall s\in{\mathbb R},$$
 where $\gamma_p:{\mathbb R}\rightarrow M$ is the integral (geodesic) curve 
of $X$ with initial condition $\gamma_p(0)=p$. In particular, $P$ must be connected.
Moreover, if $Y\in\Gamma(TP)$ is a vector field tangent to $P$, we have $Y\cdot V=
g(\nabla V,Y)=0$, because the gradient $\nabla V$ is orthogonal to the level hypersurfaces.
This means that $V({\mathcal F}(p,s))$  depends only on $s$. On the other hand, equation
(\ref{Obata}) implies that$$
(V\circ \gamma_p)'' (s) =V(\gamma_p(s)),\qquad V(\gamma_p(s))^2-\big( (V\circ\gamma_p) '(s)
\big)^2-\varepsilon^2=0,$$ 
for each $p\in P$ and $s\in{\mathbb R}$. It follows that$$
V(\gamma_p(s))=a\cosh s+b\sinh s, \quad\hbox{$b\in{\mathbb R}$\quad with }\quad
a^2=b^2+\varepsilon^2.$$
Since we are assuming that $V$ has no critical points on $M$ and $V(\gamma_p(s))$ 
only depends on $s$, we deduce that $$
(V\circ\gamma_p)' (s)=a\sinh s+b\cosh s\ne 0, \qquad \forall s\in{\mathbb R},$$
and this is equivalent to the inequality $|b|\ge a$. It turns out that $|b|=a>0$ and
$\varepsilon=0$ and so, reversing the parameter $s$ if necessary, we obtain$$
V(\gamma_p(s))=ae^s,\qquad \forall p\in P,\,\forall s\in{\mathbb R}.$$
But, from (\ref{u-r}), we know that $V(\gamma_p(s))\rightarrow +\infty$
when $\gamma_p(s)$ approaches $\partial M$. Hence $$
\exists \lim_{s\rightarrow -\infty}\gamma_p(s)=\gamma_p(-\infty)\in M\quad\hbox{and}
\quad V(\gamma_p(-\infty))=0,$$
which is a contradiction, since $V$ is positive. 

\qed 

\begin{remark}\label{r0}{\rm
If  the proof of Theorem \ref{rigidity} is closely analyzed, one can see that we
actually show that, for each Killing or parallel spinor on the spatial infinity
$\partial M$, 
there exists an imaginary Killing spinor on the  
bulk spatial manifold $(M,g)$. In fact, we have proved that there exists a parallel spinor $\Psi\in
\Gamma({\mathbb S} {\overline M})$ on the conformal compactification $({\overline M}, g^*)$.
Now, the same construction which allows to pass from parallel Euclidean spinors to
imaginary hyperbolic Killing spinors by using the conformal factor between the Euclidean and
the hyperbolic metrics on the disc $B^{n+1}$ (see \cite{BFGK,Gi}) can be used 
in order to build from $\Psi$ an imaginary Killing spinor on $(M,g)$. This means that
we can extend each supersymmetric infinitesimal isometry of $(\partial M,\gamma)$
to a supersymmetric infinitesimal isometry of the spatial slice $(M,g)$ of the
ALAdS static vacuum. This is an analogue to the
result of extending conformal transformations on the conformal infinity of a PE manifold 
to isometries on the bulk manifold proved
by M. Anderson (see \cite[Theorem 3.3]{An1} and \cite[Theorem 2.5]{An3}). In fact,
our result gives rise to the following questions: {\em Is it possible to extend each Killing vector field
on $(\partial M,\gamma)$ to a Killing vector field on $(M,g)$ for the spatial slice
of any ALAdS static vacuum? What about isometries?}
Notice that  in this supersymmetric setup for ALAdS static vacua, 
the conclusion from the existence of only one Killing or parallel spinor is
stronger than the existence of a conformal transformation in M. Anderson's result: 
the presence of only one such spinor on the boundary at infinity 
forces the bulk manifold to be maximally symmetric.  We will point out that this occurs 
because we prevent  our vacua $(M,g,V)$ to have singularities.
By the way, from the existence of the aforementioned imaginary Killing 
spinor on $(M,g)$ and the Baum's classification given in \cite{Bm} we could finish the proof
in an alternative way,
but other considerations  led us to adopt our approach.
}\end{remark}

Due to its simply-connectedness, the sphere ${\mathbb S}^{n-1}$ has a unique spin structure. Moreover, the spinor bundle corresponding to this structure and to the round Riemannian 
metric has a $2^{\left[\frac{n-1}{2}\right]}$-dimensional space of real Killing spinors for $\lambda=\frac{n-1}{2}$
and another one with the same dimension for $\lambda=-\frac{n-1}{2}$ (see \cite[p.\! 37]{BFGK},
\cite[Examples A.1.3.2]{Gi}). In this way we see that Theorem \ref{rigidity} 
 implies the rigidity result by X. Wang and by P. Chru{\'s}ciel and M. Herzlich.\\

\begin{corollary}\label{sphere}
 Let $({\mathcal V},g_{\mathcal V})$ be
an $(n+1)$-dimensional, $n\ge 3$, AAdS static vacuum spacetime determined by a static vacuum triple $(M,g,V)$. Suppose
that $M$ is a spin manifold and that the spatial slice $(\partial M,\gamma)$ of its conformal null infinity ${\mathcal J}=(\R\times\partial M,-dt^2+\gamma)$, with $\gamma= V^{-2}g_{|\partial M}$ is a unit round hypersphere. Then  $({\mathcal V},g_{\mathcal V})$ is the
AdS spacetime.\\
\end{corollary}

Besides this uniqueness of the hyperbolic space when the prescribed conformal
infinity of the spatial slices is a round sphere,  Theorem \ref{rigidity} provides a non-existence result
when this  conformal  infinity is a non-spherical compact spin manifold carrying non-trivial
real Killing spinors. In the
simply-connected case, C. B{\"a}r determined in \cite{Ba1} all these spin manifolds (see also
\cite{A} and \cite[Theorem 3]{Ba2} for non-simply-connected space forms). Using this classification, we obtain
the following: \\

\begin{corollary}\label{Killing}  
There is no $(n+1)$-dimensional, $n\ge 3$, ALAdS static vacua $(M,g,V)$ whose conformal null infinities ${\mathcal J}$ have spatial
slices $(\partial M,V^{-2}g_{|\partial M})$ isometric to  non-spherical compact spin manifolds admitting a non-trivial real Killing 
spinor. In the simply-connected case, they are  Einstein-Sasaki manifolds, 
3-Sasaki manifolds, nearly-K{\"a}hler non-K{\"a}hler 6-manifolds and 7-manifolds
carrying nice $3$-forms.
In the non-simply-connected case, they include, for example, all the round quotients 
${\mathbb S}^3/\Gamma$, where $\Gamma\subset{\mathbb S}^3$ is any of its finite subgroups,
and real  projective spaces ${\mathbb R}P^n$ with dimensions $n=8k+3$ or $n=8k+7$, $k\ge 0$.\\ 
\end{corollary}

Theorem \ref{rigidity} also provides a non-existence result
when $(\partial M,V^{-2}g_{|\partial M})$ is isometric 
to a compact spin manifold admitting non-trivial
parallel spinors.  Taking into account that the product of spin manifolds
with non-trivial parallel spinors is also another spin manifold of this same type, and that 
M. Wang determined in \cite{W1,W2} all irreducible  spin manifolds carrying 
parallel spinors (see also
\cite[Theorem A.4.2]{Gi}),  we obtain the following: \\

\begin{corollary}\label{parallel}  
There is no $(n+1)$-dimensional, $n\ge 3$, ALAdS static vacua $(M,g,V)$ whose conformal null infinities ${\mathcal J}$ have spatial
slices $(\partial M,V^{-2}g_{|\partial M})$ isometric to compact spin manifolds admitting a non-trivial parallel
spinor. In the simply-connected case, they are just Calabi-Yau manifolds, 
hyper-K{\"a}hler manifolds,  $G_2$ 7-manifolds, $Spin(7)$ 8-manifolds and
all their Riemannian products.
In the non-simply-con\-nec\-ted case, they include, for example, all the flat tori
${\mathbb T}^n$ with the trivial spin structure and all the Riemannian products 
of trivial flat tori ${\mathbb T}^k$, $1\le k\le n-2$, with the examples 
of simply-connected manifolds above. 
\end{corollary}

\begin{remark}{\rm
It is important to note that the role of spin structures is essential. 
In fact, it is well-known
that the so-called family of AdS toroidal black hole metrics (see \cite{BMN} or \cite[Example 2.2, Remark 3.4 ii)]
{An1}), constructed on the solid $n$-dimensional torus $B^2\times {\mathbb T}^{n-2}$  are  the unique  PE 
manifolds whose conformal infinity is the flat torus ${\mathbb T}^{n-1}$. They are given by$$
g=\frac1{U(r)}dr^2+U(r)d\theta^2+r^2g_{{\mathbb T}^{n-2}},$$
where $g_{{\mathbb T}^{n-2}}$ is the standard flat metric on the torus and$$
U(r)=r^2\left(1-\frac{r_0^n}{r^n}\right),$$ 
with $r_0>0$ and $\theta$ a $(4\pi/nr_0)$-periodic coordinate. One can check that,
taking $V=V(r)=r^2$, we have that $(B^2\times{\mathbb T}^{n-2},g,V)$ is
a static vacuum. Indeed it determines the so-called AdS soliton (see \cite{Wa2}
and \cite{GSW} for a uniqueness result).
It is clear that $M=B^2\times{\mathbb T}^{n-2}$ is a spin manifold, but the spin structure inherited
on its boundary at infinity $\partial M={\mathbb T}^{n-1}$ is not the trivial one, since $M$ is constructed by gluing on a $2$-disc
onto a simple closed geodesic of the flat ${\mathbb T}^{n-1}$. Then, even though $(\partial M,\gamma)$ is a 
flat torus, the spin structure inherited from $(M,g)$ admits no parallel spinor fields and Corollary \ref{parallel}
does not apply in this situation. \\ 
}\end{remark}

\section{Supersymmetries on the boundary and spin conformal  compactifications with singularities}
If we compare Theorem 3.3 in \cite{An1} (see also Theorem 2.5 in \cite{An3}) with the main Theorem 
\ref{rigidity} of this paper, we see that, in the spirit of the AdS/CFT correspondence, 
all the conformal transformations on the
 boundary of a PE manifold come from 
Riemannian isometries of its conformal compactifications and that, however, in the case
of ALAdS static vacua, the supersymmetric infinitesimal isometries of the null infinity
give non-existence results, except for the spherical case.\\

Hence, the proof of Theorem \ref{rigidity} can be read as a supersymmetric version of
Anderson's result under different conditions on the Ricci curvature. Indeed, its proof shows that each Killing
or parallel spinor on the spatial slice of the null infinity
comes from an imaginary Killing spinor on the bulk manifold (see Remark \ref{r0}). But a careful reading of this proof 
makes also manifest that, if we  exclude conical or cusp 
singularities, the only vacua supporting these imaginary Killing spinors are
AdS spacetimes. This is why we finish this paper with some examples showing that
there are ALAdS static vacua with null infinities carrying Killing or parallel spinors  provided
that  singularities are allowed. \\

We could say that  supersymmetries on
a non-spherical null infinity yield static vacua with hyperbolic
conical or cusps singularities.
Mathematically, these examples take the form of warped products. From a
physical point of view such kind of manifolds also appear regularly (\cite{AC,BMN}). On the other hand, a strong necessity of
considering ALAdS metrics with singularities arises as well when
one wants to understand the manifold structure of the space of these metrics, with a given
topology (see, for example, \cite[(3.6)]{An1}).\\ 

{\bf Example 1.} Consider  the Poincar{\'e} hyperbolic ball $\left(
B^{n},\frac{4|dx|^2}{(1-|x|^2)^2}\right)$ as a model of the hyperbolic space. Using polar coordinates $x=rp$, with
$r\in ]0,1]$ and $p\in{\mathbb S}^{n-1}$, we have that the hyperbolic metric of constant 
sectional curvature $-1$ takes
the form$$
g=\left(\frac{2}{1-r^2}\right)^2\left(dr^2+r^2\gamma_{{\mathbb S}^{n-1}}\right),$$
where $\gamma_{{\mathbb S}^{n-1}}$ stands for the unit round metric on the sphere. 
If we consider the change of variables given by $s=\ln\frac{1+r}{1-r}\in{\mathbb R}^+$, we
obtain$$
g=ds^2+(\sinh^2 s)\gamma_{{\mathbb S}^{n-1}}.$$
These two expressions for the Poincar{\'e} metric are valid only on the punctured 
ball $B^{n}-\{0\}\cong\,]0,1[\times {\mathbb S}^{n-1}\cong {\mathbb R}^+\times {\mathbb S}^{n-1}$, although
they are smoothly extendable to the origin. This latter is an example of the so-called
{\em warped Riemannian products} (see, for instance, \cite{Be,O'N,K}). \\

In general, if $I\subset {\mathbb R}$ is an open 
interval, $(P,\gamma)$ a Riemannian $(n-1)$-manifold and $f\in C^\infty(I)$ is a 
positive function, we will say that the $n$-dimensional Riemannian manifold
$(I\times P, g=ds^2+f(s)^2\gamma)$ is the  product of $I$ and $P$  warped by means of
the function $f$. For the sake of simplicity, if we only consider  Einstein manifolds, using the form of the Ricci tensor on a
warped product (see \cite[Lemma 4]{K}) and recalling (\ref{scalar}), we have restrict ourselves to warping functions
$f$ satisfying the linear differential equation $f''-f=0$. With this choice, we ensure  that $\Ric_g(\frac{\partial}{\partial s},
\frac{\partial}{\partial s})=-(n-1)$ at each point of $I\times P$. Taking now into account 
the values of $\Ric_g$ on the directions orthogonal to the vector field $\frac{\partial}
{\partial s}$, that is, directions tangent to $P$, we conclude that there are essentially three types of warped products 
which eventually may produce Einstein manifolds with scalar curvature $-n(n-1)$. The
first one is the {\em hyperbolic cone} on a given compact Riemannian manifold $(P,\gamma)$, given by$$
({\mathbb R}^+\times P,g=ds^2+(\sinh^2 s)\gamma).$$ 
It is immediate to see again from \cite[Lemma 4]{K}, for instance, that we have for the 
directions tangent to this cone and perpendicular
to the radial direction $$
\Ric_g+(n-1)g=\Ric_\gamma-(n-2)\gamma.$$ 
So, we will assume that $(P,\gamma)$ is an Einstein manifold with scalar curvature $(n-1)(n-2)$. 
We can also see in \cite{K} that the smooth function defined on $\R^+\times P$ by $V(s,p)=
\cosh s$, for each $s\in{\mathbb R}^+$ and $p\in P$, is a solution to the Obata type equation $\nabla^2V=Vg$. 
Thus the triple $(M,g,V)=(\R^+\times P,ds^2+(\sinh^2r)\gamma,\cosh r)$ determines a static vacuum. Defining a 
new variable $t\in ]0,\frac{\pi}{2}]$ by the equality  $t=-\frac{\pi}{2}+2\arctan e^s$,  we obtain
$$
{\overline g}=\frac1{\cosh^2 s}g=dt^2+(\sin t^2)\gamma.$$
The spherical conical metric $dt^2+(\sin t^2)\gamma$ obviously extends to $[0,\frac{\pi}{2}]\times P$, a compact
 manifold with boundary $\{\frac{\pi}{2}\}\times P\cong P$ and a
conical singularity at $t=0$. This singularity is removable if and only if $(P,\gamma)$ is the round unit $(n-1)$-sphere
and, in this case, the corresponding hyperbolic cone is nothing but the $n$-dimensional hyperbolic space
(see \cite[p. 269, Lemma 9.114]{Be}). When $P$ is chosen to be one of the non-spherical compact spin $(n-1)$-dimensional 
manifolds listed in Corollary \ref{Killing}, we get an example of non-AdS ALAdS static vacuum with spatial
infinity supporting non-trivial Killing spinors with a conical singularity.\\

{\bf Example 2.} The second type of warped product relevant to our purposes is the so-called {\em 
hyperbolic cusp} on a compact Riemannian manifold $(P,\gamma)$, given by$$
({\mathbb R}\times P,g=ds^2+e^{2s}\gamma).$$
For these cusps, the Ricci curvature is also $-(n-1)$ along the radial direction $\frac{\partial}
{\partial s}$ and, for orthogonal directions tangent to $P$, we have$$  
\Ric_g+(n-1)g=\Ric_\gamma.$$
Then, we will assume that in this case $(P,\gamma)$ is Ricci-flat. As in Example 1, we can see in \cite{K} that $V(s,p)=e^s$ is a solution to $\nabla^2V =Vg$. Defining a new variable $t\in{\mathbb R}^+$ as $t=e^s$, we obtain$$
{\overline g}=\frac1{e^{2s}}g=dt^2+\gamma.$$
The cylindrical metric $dt^2+\gamma$ clearly extends to $[0,+\infty[\times P$, which is a non-compact
 manifold with boundary $\{0\}\times P\cong P$ and one cylindrical end. So, the corresponding
 static vacuum $(M,g,V)=(\R^+\times P,ds^2+e^{2s}\gamma,e^s)$ has a null infinity whose
 Riemannian slice is $(P,\gamma)$  for
$s=+\infty$, but it is not compactifiable at $s=-\infty$ because of the presence of a hyperbolic cusp. 
Of course, the cusp singularity is
always unremovable, although it has finite volume. When $P$ is chosen to be one of the non-spherical compact spin $(n-1)$-dimensional 
manifolds listed in Corollary \ref{parallel}, we get an example of non-AdS ALAdS static vacuum with spatial
infinity supporting non-trivial parallel spinors having a cusp singularity.\\

\end{document}